\documentclass[a4paper,12pt]{article}
\usepackage{amsmath,amsfonts,amssymb}
\usepackage{a4wide}
\usepackage{hyperref} 
\usepackage{mathrsfs}

\usepackage{fancyhdr} 
\makeindex 

\newcounter{EQNR}[N]
\newcommand{\nne}{\refstepcounter{EQNR} \tag{\theEQNR}}

\newcommand{\SO}{{\rm SO}}

\begin{document}
$$\textbf{\Large Round about Kudla's Green function for $\textrm{SO}(3,2)$}\\$$

$$\textbf{Rolf Berndt}$$\\
\begin{abstract} 
Kudla conjectured that certain Eisenstein series 
contain important arithmetical and geometric information. 
The following note describes a certain aspect of this general picture in the special case concerning the orthogonal group $\textrm{SO}(3,2).$ \\
\end{abstract} 

This is a short version of a long report 'Around Kudla's Green function for $\textrm{SO}(3,2)'$, which contains much more details and 
can be found on the authors home page\footnote{\url{https://www.math.uni-hamburg.de/en/forschung/bereiche/ad/personen/berndt-rolf.html} }.
The main result is the explicit determination of the integral of Kudla's Green function for $\textrm{SO}(3,2)$. 
In the long file there are also results concerning $\textrm{SO}(p,2)$ for $p=2 \,\textrm{and}\,1$. All this work was insprired by my wish to  extend  the results for  $\textrm{SO}(p,2)$ obtained in joint work work with Ulf K\"uhn  \cite{BeKI}, \cite{BeKII},\cite{BeKIII}. 
\\

{\textbf{A dual pair.}}\\

We start by looking at the 
 (special) \textit{dual} pair of groups $G$ and $H,$ where 
$G = \textrm{SL}(2,\mathbb{R)},$ 
the group of two by two real matrices $g = \begin{pmatrix}a&b\\c&d \end{pmatrix}$ with determinant one, 
and $H$ an orthogonal group stabilizing a lattice $L$ 
or a quadratic form given like this: 
We take an $\mathbb{R}-$vector space $V$ of dimension $n$ 
equipped with a bilinear form $(x,y) = \sum_{i,j=1}^n 2a_{ij}x_iy_j = {}^txAy,$ 
 quadratic form $q(x) = (1/2)(x,x)$ and (special) lattice $L \simeq \mathbb{Z}^n$ 
 of signature $(p,q), p+q=n$ with discriminant $d = \textrm{det}\,A;$ for instance  
$H = \textrm{SO}(p,q)$ the group of $n$ by $n$ real matrices $M$ with determinant 
one and ${}^tME_{pq}M = E_{pq} = \begin{pmatrix}E_p&\\&-E_q \end{pmatrix}$\\

$G$ acts on the upper half plane $\mathbb{H} = G/\textrm{SO}(2) = \{\tau = u+iv \in \mathbb{C};v>0\}$
$$
G \times \mathbb{H} \rightarrow \mathbb{H}, \quad (g,\tau) \mapsto g(\tau) = \frac{a\tau +b}{c\tau +d}.
$$
$H$ with maximal compact subgroup $K$ has as homogeneous space $\mathbb{D} = H/K$ 
where one uses many different realizations, for instance, in the case of $p=2$ or $q = 2$ which will be the most 
important for us, the Grassmannian $\textrm{Gr}_2(V)$ of two-dimensional oriented negative 
planes $X \subset V.$\\

\pagebreak
 \textbf{Eisenstein series}\\

There are different approaches to generalize the classical Eisenstein series 
mainly using the adelic language. But here we stay in the archimedean world 
and following the paper \cite{BK} by Bruinier and K\"uhn present the following generalization: 
Let $L \subset V$ be an even lattice, $L'$ its dual, and $\Gamma(L)$ the discriminant group, 
in particular:\\
\textbf{Example (2,3).} $V = \mathbb{R}^5$ equipped with the quadratic form of signature (2,3)
$$
q(\textbf{x}) = x_1x_5 + x_2x_4 - x_3^2,
$$
and the lattice $L = \mathbb{Z}^5,$ hence $L' = \{(y_1,y_2,y_3,y_4,y_5); y_i\in \mathbb{Z}, i\not=3, y_3\in (1/2)\mathbb{Z}\},$ and $\Gamma(L) \simeq \mathbb{Z}/2\mathbb{Z}
\simeq \{\beta = \bar{0},\bar{1} \}$\\
Using the \textit{Weil representations} of $\textrm{SL}$ resp. its metaplectic cover $\textrm{Mp}$ 
in Section 3 of \cite{BK} (and Section 5 of [Ku1]) one finds vector valued 
Eisenstein series $E_\beta$ of weight $\ell \geqslant 2, \ell \in (1/2)\mathbb{Z}$ (indexed by the 
elements $\beta$ of $L'/L$) with 
Fourier expansion
\begin{align*}
E_\beta(\tau,s) 
&= \sum_{\gamma\in L'/L}\sum_{m\in \mathbb{Z}-q(\gamma)}c_\beta(\gamma,m,s,v) e^{2\pi imu}e_\gamma\\
&= \sum_{\gamma\in L'/L}\sum_{m\in \mathbb{Z}-q(\gamma)}C_\beta(\gamma,m,s)\mathcal{W}_s(a) e^{2\pi imu}e_\gamma\, \quad a = 4\pi mv.
\end{align*}
For the case of the example (2.3) above, from [BK] Theorem 3.3 we have  formulae for $c_\beta(\gamma,m,s,v)$ 
which we don't reproduce here, but only for $\beta = \bar{0} = 0$ and
\begin{align*}
\label{Sint1}\nne
c_0(\gamma,m,0,v) 
&= \begin{cases} 
 C(\gamma,m,0) e^{-a/2} \,\textrm{for}\,m>0,\, \quad a = 4\pi mv,\\
  0, \quad\,\textrm{for}\,m<0,
\end{cases} \\
 C(\gamma,m,0):&= - 2^6\cdot 3\cdot 5 \cdot \pi^{-2} |m|^{3/2}L(\chi_{D_0},2)\sigma_{\gamma,m}(5/2). 
\end{align*}
Here $L(\chi,s) = \sum_{n=1}^\infty \frac{\chi(n)}{n^s}$ is the usual 
$L-$function here for the character $\chi = \chi_{D_0},$ where $m$ and $D_0$ are related by 
\begin{align*}
\label{mA}\nne
D_0 f^2 &= 4m \,\textrm{for}\,m\in \mathbb{Z}\,\,\textrm{i.e.,}\,\textbf{Case A}\\
&= 4^2m \,\textrm{for}\,m-1/4 \in \mathbb{Z}\,\,\textrm{i.e.,}\,\textbf{Case B},
\end{align*}
as in \cite{BK} (3.24).
And $\sigma_{\gamma,m}(5/2)$ is the  generalized divisor sum 
from \cite{BK} (3.28).\\

We abbreviate the derivative $\frac{d}{ds}$ by the asterisque $'$ 
and by some nontrivial manipulations get the final formulae
\begin{align*}
\label{SGint2}
c'_0(\gamma,m,0,v) 
&=
\begin{cases}
 C(\gamma,m,0) e^{-a/2}(J_+(3/2,a) + \frac{C'(\gamma,m,0)}{C(\gamma,m,0)} ) ,\,\textrm{for}\,m>0,\\
 C(\gamma,m,0) e^{-|a|/2} \cdot J_-(3/2,|a|) ,\,\textrm{for}\,m<0,
\nne
\end{cases}
\end{align*}
with
\begin{align*}
\label{FEisl08}
J_+(s,a) &:= \int_0^\infty e^{-aw}((w+1)^s - 1)dw/w,\\
J_-(s,a) &:= \int_0^\infty e^{-aw}w^s dw/(w+1).
\nne
\end{align*}
In the following, we try to describe some arithmetic and geometric content contained in these formulae.
Up to now the orthogonal group $H$ didn't appear much:
 in the special case $H \simeq \textrm{SO}(3,2)$ described above, it entered via the two dimensional 
 Weil representation of the metaplectic group. It shall become much more present now.\\
 
 \textbf{Our orthogonal world.}\\
 
  There is a second way to the formulae for the Eisenstein Fourier coefficients  which has a geometric touch,
 As said above, for $q=2,$ the symmetric space associated with 
 our orthogonal group $H$ may be identified with the
set 
\begin{align}
\label{Dp2}\nne 
\mathbb D=\{\mbox{ oriented negative 2-planes } X \subset V = \mathbb R^n\,\},
\end{align}
i.e. $ X = \,<v_1,v_2>,\, v_j \in V, \,q(v_j) < 0, \,(v_1,v_2) = 0.$  
 It is well-known that, for $p=3,$ $\mathbb D$ has
two connected components $\mathbb D^+$ and $\mathbb D^-$, and $\mathbb D^+$ is isomorphic to
the Siegel half plane  $\mathbb H_2$ of genus 2, for $p=2$ one has $\mathbb D ^+\simeq \mathbb H\times \mathbb H,$ and $\mathbb D ^+\simeq \mathbb H$ for $p=1.$ 
These facts have their source in special \textit{exceptional} homomorphisms 
between the groups 
\begin{align*}
\label{Dp25}\nne 
\textrm{SL}(2,\mathbb{R}) \rightarrow \textrm{SO}(2,1)\\
\textrm{SL}(2,\mathbb{R}) \times \textrm{SL}(2,\mathbb{R})\rightarrow \textrm{SO}(2,2)\\
\textrm{Sp}(2,\mathbb{R}) \rightarrow \textrm{SO}(3,2)
\end{align*}
The isomorphism $\mathbb{H}_2 \rightarrow \mathbb{D}^+$ says 
that $\mathbb{D}^+$ may be parametrized by three complex 
variables $z_1,z_2,z_3$ with $y_1y_3 - y_2^2 > 0.$ 
In our parametrization by these $z = \begin{pmatrix}z_1&z_2\\z_2&z_3 \end{pmatrix} \in \mathbb{H}_2,$ elements of $\mathbb{D}$ are negative planes 
$X_z = <\textrm{Re}\,u(z), \,\textrm{Im}\,u(z)>, z \in \mathbb{H}_2$ with 
\begin{align}
\label{Zxx} \nne
u_1(z) = z_2^2-z_1z_3, u_2(z) = -z_1, u_3(z) = -z_2, u_4(z) = -z_3, u_5(z) = 1.
\end{align}

Starting from $\mathbb{D} = \textrm{Gr}_2(V),$ we come to \textit{quasiprojective algebraic varieties} $X_\Gamma = \Gamma\backslash \mathbb{D}$ 
made by the $\Gamma-$orbits, where $\Gamma(L) = \textrm{SO}(L)$ is the discrete subgroup 
preserving the lattice $L.$
Given a variety or any space, one is tempted to look for curves or subspaces. 
In our case, the first natural objects are \textit{divisors} coming up as follows. 
For any vector $x \in V$ with $q(x) = m >0$ 
the orthogonal complement in $\textrm{Gr}_2(V)$ defines a divisor $Z(x)$
on $\textrm{Gr}_2(V).$ In the coordinates above
\begin{align*}
\label{Zxxx} \nne
Z(x) &= \{z \in \mathbb{D}(V); z \bot x \},\\
&\simeq \{z \in \mathbb{H}_2; {}^tu(z)\tilde{Q}x = 0 \},\\
&\simeq \{z \in \mathbb{H}_2; \psi _{u(z),x} = x_1-x_2z_3+2x_3z_2-x_4z_1+x_5(z_2^2-z_1z_3) = 0 \}.
\end{align*}
This is a two-dimensional complex subvariety in $\mathbb{H}_2$ called 
a \textit{Humbert surface} (of discriminant $\Delta = (2x_3)^2 - 4x_1x_5 - 4x_2x_4$).
If $\beta \in L'/L$ and $m \in \mathbb{Z} + q(\beta) > 0,$ 
then 
$$
\mathcal{H}(\beta,m) = \sum_{x\in L+\beta, q(x)=m} x^\bot
$$
is a $\Gamma(L)-$invariant divisor on $\textrm{Gr}_2(V),$ called the 
\textit{Heegner divisor} of discriminant $(\beta.m).$ It gets the same name as divisor on $X_\Gamma.$ 
For a divisor $D$ on $X_\Gamma$ with volume form $\omega,$ 
one defines
$$
\textrm{deg}(D) = \int_D \omega.
$$
In \cite{BK} Proposition 4.8 or \cite{Ku1} Proposition 5.1 by some highly nontrivial manipulations 
the first mysterious and wonderful result was obtained: \\

\textbf{Theorem 1.} For the (2,3) case from the example above, one has
\begin{align*}
\label{epilo1}\nne
\frac{2}{B}\textrm{deg}\,(\mathcal{H}(\gamma,m)) &= - C_0(\gamma,m,0).\\
\end{align*}
Here, we have (a well known result going back to Siegel)
\begin{align*}
\label{epilo2}\nne
B = \int_{X_L} \Omega^3 = \zeta(-1) \zeta(-3) =  (1/12)(1/120) = 2^{-5}\cdot 3^{-2}\cdot 5^{-1}
\end{align*}
with ([Ku1] (5.8))
\begin{align*}
\label{epilo3}\nne
\Omega^3 = - \frac{3}{16\pi^3}(y_1y_3-y_2^2)^{-3}(\frac{i}{2})^3
dz_1\wedge d\bar{z}_1\wedge dz_2\wedge d\bar{z}_2\wedge dz_3\wedge d\bar{z}_3.
\end{align*}
I.e., we have from (\ref{Sint1})
\begin{align*}
\label{SGint1a}\nne
\textrm{deg}\,(\mathcal{H}(\gamma,m)) =  - (B/2) C(\gamma,m,0) =
 -  3^{-1} \cdot \pi^{-2} |m|^{3/2}L(\chi_{d_F},2)\sigma_{\gamma,m}(5/2). 
\end{align*}
 {\bf  Remark 1.} \label{zuu3} There is another way to see the result of Theorem 1:
 
  To compare, we also reproduce Kudla's treatment. In \cite{Ku1} p.337ff, things look like this. 
We have a two-component Eisenstein series
\begin{align*}
\label{kud1}\nne
E(\tau,3/2;\varphi_0) &= 1 + \zeta(-3)^{-1}\sum_{m=1}^\infty H(2,4m)q^m = \sum A_0(m,v)q^m\\
E(\tau,3/2;\varphi_1) &=  \zeta(-3)^{-1}\sum_{m-1/4=0}^\infty H(2,4m)q^m = \sum A_1(m,v)q^m\\
\zeta(-3)^{-1} &= 2^3\cdot3\cdot5,\quad 4m = D_0f^2, D_0 \equiv  0, 1 \,\textrm{mod}\,4 \quad\,(\textrm{\cite{Ku1}}(5.18)),\\
  H(2,4m) &= 
L(-1,\chi_{D_0})\sum_{d|f}\mu(d)\chi_{D_0}(d)d\,\sigma_3(f/d)\\ 
 &=  - (1/(2\pi^2))L(2,\chi_{D_0})D_0^{3/2}\xi(D_0,f),\,
\xi(D_0,f) = \sum_{d|f}\mu(d)\chi_{D_0}(d)d\,\sigma_3(f/d).
\end{align*}
Hence, Kudla's Eisenstein series is half of the series in Bruinier-K\"uhn and here we have the coefficients 
\begin{align*}
\label{kud2}\nne
A(m,v) &= 2^3\cdot3\cdot5\cdot
H(2,4m)
 =  2^6\cdot3\cdot5\cdot
  L(2,\chi_{D_0})m^{3/2}/(2\pi^2)f^{-3}\xi(D_0,f).
\end{align*}
Finally from \cite{Ku1} (5.22), one has 
\begin{align*}
\label{ku4}\nne
\textrm{deg}\,(Z(m,\varphi_\mu)) &= \textrm{deg}\,G_{4m} = - (1/12) H(2,4m)\\
&= (1/(24\pi^2)) \cdot D_0^{3/2}L(\chi_{D_0},2)\xi(D_0,f).
\end{align*}
 
 \textbf{Kudla's Green function.}\\
 
There is still another even more surprising result interpreting the Fourier coefficients $c'(\gamma,m,0,v)$ 
of the derived Eisenstein series via certain geometric integrals:\\

The classical notion of a Green function has been refined in several ways.
For us, given the space $X_\Gamma$ and a divisor $D$ on $X_\Gamma$ a Green 
function $\Xi$ is real analytic on $X_\Gamma \backslash D$ and has a logarithmic 
singularity along $D.$ We will describe \textit{Kudla's Green function} which uses the \textit{majorant} 
of a quadratic form: Given a symmetric $Q \in \textrm{M}_n(\mathbb{Q})$ and a quadratic form ${}^txQx =: Q[x]$ (Siegel's notion), a majorant is given by a positive definite symmetric $P\in \textrm{M}_n(\mathbb{Q})$ 
minimal with ${}^txPx \geqslant {}^txQx$ for all $x.$ On the first pages of \cite{S1}, Siegel 
shows that $P$ defines such a majorant if
\begin{align}
\label{mj1}
 \nne
 PQ^{-1}P = Q, \, P = {}^tP > 0.
\end{align}
If  $G\subset \textrm{GL}_n(\mathbb{R})$ is the group preserving $Q$ 
with maximal compact $K,$ these $P$ are parametrized by $\mathbb{D} = G/K.$ 
In Kudla's notation, for $z \in \mathbb{D}$ and $(x,x) = {}^txQx,$ one writes 
for the majorant in $z$ 
\begin{align}
\label{mj2}
 \nne
 (x,x)_z = (x,x) +2R(x,z), R(x,z) = - (\textrm{pr}_z(x),\textrm{pr}_z(x)),
\end{align}
where $\textrm{pr}_z: V = \mathbb{R}^n \rightarrow z \in \mathbb{D} = \textrm{Gr}_2(V)$ 
the projection with kernel $z^\bot.$ For our example (3,2) $q(x) = x_3^2 - x_1x_5 - x_2x_4 = (1/2)(x,x)$ 
we come to 
\begin{align}
\label{mj3}
 \nne
 R(x,z) = R(x,X_z) = \frac{|x_1 - x_2z_3 + 2x_3z_2 - x_4z_1 + x_5(z_3 - z_1z_2)|^2}{2(y_1y_3 - y_2^2)} 
 = \frac{|\psi_{u(z),x}|^2}{2\eta^2}.
\end{align}
In \cite{Ku1} (4.28), Kudla introduces
\begin{align}
\label{mj4}
 \nne
 \xi(x,z) := \beta_1(2\pi R(x,z)) e^{-2\pi Q(x)},\,\,\beta_s(x) = \int_1^\infty e^{-xt}dt/t^s
\end{align}
which has a logarithmic singularity as $t\rightarrow 0$ and decays like $e^{-t}$ as $t\rightarrow \infty.$ 
In \cite{Ya} (4.5), Yang takes $\beta = \beta_1$ and a lattice $L$ to 
construct the Green function 
\begin{align}
\label{mj5}
 \nne
 \Xi(m,v) := \frac{1}{2}\sum _{0\not=x\in L(m)}\beta(x\sqrt v,z), m\in \mathbb{Z}, L(m) = \{x \in L; q(x) = m\}.
\end{align}
We stick to this definition but have to modify the lattice-summation to adapt to our situation.
Above, we had $X_z =\, \langle\textrm{Re}\,u(z), \,\textrm{Im}\,u(z)\rangle, z \in \mathbb{H}_2$ with (\ref{Zxx})
\begin{align*}
u_1(z) = z_2^2-z_1z_3, u_2(z) = -z_1, u_3(z) = -z_2, u_4(z) = -z_3, u_5(z) = 1.
\end{align*}
Hence, 
\begin{align*}
Z(x) &= \{z \in \mathbb{D}; z \bot x \},\\
&\simeq \{z \in \mathbb{H}_2; {}^tu(z)\tilde{Q}x = 0 \},\\
&\simeq \{z \in \mathbb{H}_2; \psi _{u(z),x} = x_1-x_2z_3+2x_3z_2-x_4z_1+x_5(z_2^2-z_1z_3) = 0 \}.
\end{align*}

 Guided by Gritsenko-Nikulin \cite{GN} p.204 and \cite{Ku1} p.338, on $\mathcal{X} := \Gamma \backslash \mathbb{D},$ we find the divisors $Z(0.m)$ and $Z(1,m)$ induced by
\begin{align}
\label{Zxi} \nne
 \sum_{x\in 2L', \tilde{q}(x) = 4m} Z(x) \quad\,\textrm{resp.}
 \sum_{x\in 2L', \tilde{q}(x) = 4M+1} Z(x) 
 \end{align}
 with $m,M \in \mathbb{Z}$ resp. $m = M+1/4$ in the second case,
which have the Humbert equation
\begin{align}
\label{Zxii} \nne
\psi _{u(z),2x} = 2(x_1-x_2z_3+x_3z_2-x_4z_1+x_5(z_2^2-z_1z_3)) = 0
 \end{align}\textsl{}
We abbreviate $L_m := \{x\in 2L', \tilde{q}(x) = 4m\}.$ 
From the discussion of the Fourier coefficients of the Eisenstein series we know 
that we have the two cases $m \in \mathbb{Z}$ and $m = M+1/4 \in \mathbb{Z} + 1/4.$ 
To better distinguish, for the  second case, we write 
\begin{align}
\label{Zxiii} \nne
 L^1_\Delta := \{x\in 2L', \tilde{q}(x) = 4m = 4M+1\}, \, \Delta = 4m = 4M+1 \,\textrm{if}\,\, m = M+1/4
 \end{align}
We introduce a two component Kudla Green function 
\begin{align*}
\label{Xxi} 
\Xi(0,m,v,z) &= \sum_{x\in L_m}\beta(2\pi vR(z,x)),\\
\Xi(1,m,v,z) &= \sum_{x\in L^1_m}\beta(2\pi vR(z,x)).\nne
\end{align*}
\textbf{The Green function integrals}\\
 
 As a next step beyond Theorem 1, Kudla in \cite{Ku1}
proposed to look at the integral of the Green function. 
Hence, we want to compare the $s-$derivatives of the Fourier coefficients  $c(\gamma,m,s,v)$ (in $s=0$) of the Eisenstein series 
with the Green function integrals 
\begin{align*}
\label{Xx2} 
I(0,m,v) :&= \int _X \Xi(0,m,v,z) d\mu_z = \sum_{x\in L_{0,m}}\int_X\beta(2\pi vR(z,x))d\mu_z,\\
I(1,m,v) :&= \int _X \Xi(1,m,v,z) d\mu_z= \sum_{x\in L_{1,m}}\int_X\beta(2\pi vR(z,x))d\mu_z.\nne
\end{align*}
By evaluating these integrals, with (\ref{SGint2}) and (\ref{FEisl08}), we get (up to a certain factor) in conformity with one of Kudla's conjectures:\\

 \textbf{Theorem 2.} For  the (3,2) case, one has
\begin{align*}
\label{Apen80}
C(\gamma,m,0) &:= - 2^6\cdot 3\cdot 5 \cdot \pi^{-2} |m|^{3/2}L(\chi_{d_F},2)\sigma_{\gamma,m}(5/2),\\
c_0(\gamma,m,0,v) 
&=   \begin{cases}C(\gamma,m,0) e^{-a/2} \,\textrm{for}\,m>0,\\
 0, \quad\,\textrm{for}\,m<0, 
 \end{cases}\\
c'_0(\gamma,m,0,v) 
&=\begin{cases} C(\gamma,m,0) e^{-a/2}(J_+(3/2,a) + \frac{C'(\gamma,m,0)}{C(\gamma,m,0)} ) ,\,\textrm{for}\,m>0,\\
 C(\gamma,m,0) e^{-|a|/2} \cdot J_-(3/2,a) ,\,\textrm{for}\,m<0,
 \end{cases}
 \\
(4/B)\cdot I(\gamma,m,v)
&=\begin{cases}  C(\gamma,m,0)J_+(3/2,a), \,\textrm{for}\,m>0,\\
  C(\gamma,m,0)J_-(3/2,a) e^{-|a|}, \,\textrm{for}\,m<0.
  \end{cases}\\
  (4/B)\cdot I^{BK}(\gamma,m,v)
&=\begin{cases} 
 - C(\gamma,m,0)\big(\frac{C'(\gamma,m,0)}{C(\gamma,m,0)} + \textrm{log}(4\pi) - \Gamma'(1) \big), \,\textrm{for}\,m>0,\\
   0, \,\textrm{for}\,m<0.
  \end{cases}
\nne
\end{align*}
Here $I^{BK}(\beta,m,v) = \int_X G_{\beta,m}(Z) \Omega^3$ is the integral of the Green function $G_{\beta,m}(Z)$
from Bruinier-K\"uhn \cite{BK} Definition 4.5 and Theorem 4.11.\\

 \textbf{Corollary.} We have
\begin{align*}
\label{Apen81}
c'_0(\gamma,m,0,v) = e^{-a/2} ((4/B)\cdot( I(\gamma,m,v) - I^{BK}(\gamma,m,v)) + \ast\, \,c_0(\gamma,m,0,v) ).
\nne
\end{align*}

\textbf{Proof of Theorem 2.}\,(Sketch)\,
 We only reproduce some elements of our calculation.
We have
\begin{align*}
\label{epi3}\nne
L_{\gamma,m} &= \{u\in \mathbb Z^5; \hat{q}(u)= u_3^2-4u_2u_4-4u_1u_5= 4m \},\\
&= \sum_{n^2|\delta_\gamma m} nL^\ast_{\gamma,m/n^2},
L^\ast_{\gamma,m} = \{u \in L_{\gamma,m}; \textrm{gcd}(u) = 1 \}\\
\gamma & = 0\,\,\textrm{for}\,\,m \in \mathbb{Z},\,\,
\gamma  = 1 \,\,\textrm{for}\,\,m = M + 1/4 \in \mathbb{Z} + 1/4\\
\delta_\gamma &= 1 \,\textrm{or}\, 4\,\quad\textrm{for}\,\,\gamma = 0\, \,\textrm{or}\,\,1.
\end{align*}

From the general theory of Humbert surfaces, we know that $\Gamma = \textrm{Sp}(2,\mathbb{Z})$ acts transitively  on quintuples $(u_1,\dots,u_5)$ without common divisor.  We have
\begin{align*}
\label{epi4}\nne
L^\ast_{\gamma,m} &= (\Gamma/\Gamma_{a})a, \\
a &= a_{\gamma,m} = {}^t(1,0,0,0,-m) \,\,\textrm{or} \,= {}^t(0,1,1,-(4m-1)/4,0) \,\textrm{for}\,\,\gamma = 0\,\, \textrm{or}\,\,1, 
\end{align*}
where $\Gamma_{a}  $ fixes $ a. $ Using this, we infer
\begin{align*}
\label{epi1a}
I(0,m,v) &= \sum_{n^2|\delta_\gamma m}I(0,v,m,n)\\
I(0,m,v,n) &:= \int _{\Gamma_{a_{0,m/n^2}} \setminus \mathbb H_2 } \int_1^\infty
e^ {-2\pi vR(na_{0,m/n^2},z)r}dr/rd\mu (z).
\nne
\end{align*}
Similarly, for $I(1,m,v).$\\

 \label{Hepi1a} We have to determine $I(\gamma,m,n,v)$ 
and follow Kudla
in \cite{Ku1} p.318. There, he treats the more general ($p$,2)-case and $x \in V(\mathbb{R})$ 
with $Q(x) = m:$ For $m>0,$ 
one chooses a basis $\textbf{v}$ for $V(\mathbb{R})$ so that the inner product 
has matrix $I_{p,2}$ and so that the respective special element $x$
is a nonzero multiple of the first basis vector $v_1,$ i.e. $x = 2\alpha v_1.$ Then 
$\SO(V)(\mathbb{R})^+ = \SO^+(p,2) = G$ and the subgroup 
stabilizing $x$ is isomorphic to  $\SO^+(p-1,2) = G_{x}.$ 
Kudla further proposes $z_0 \in \mathbb D$ to be 
the oriented negative 2-plane spanned by $v_{p+1}$ and
$v_{p+2}$ and let $K = \SO(n)\times \SO(2)$ be its 
stabilizer in $G.$ The plane spanned by $v_1$ and $v_{p+1},$ 
the first negative basis vector, has signature (1,1). 
The identity component of the special orthogonal 
group of this plane is a 1-parameter subgroup   
\begin{align*}
\label{Gint4}
A = \{a_t =
\left(\begin{smallmatrix}\cosh t &\sinh t \\ \sinh t &\cosh t \end{smallmatrix}\right) ; t \in \mathbb{R} \}
\nne
\end{align*}
where $a_t \cdot v_1 = v_1 \cosh t  +  v_{n+1} \sinh t.$ 
Let $A_+$ the subset of $a_t'$s with $t\geq 0.$
Then, from the general theory of semisimple 
symmetric spaces - with \cite{FlJ} as a reference - Kudla has a 'double set decomposition' 
(as in Bruinier and Funke \cite{BF1} Section 7)
\begin{align*}
\label{Gint5}
G = G_{x}A_+K
\nne
\end{align*}
and the integral formula
\begin{align*}
\label{Gint6}
\int_G \Phi (g) dg = \int_{G_{x}}\int_{A_+}\int_K
\Phi (g_xa_tk)\lvert \sinh t \lvert (\cosh t)^{p-1} dg_xdtdk.
\nne
\end{align*}
For $z = g_{x}a_t \cdot z_0,$ one has 
\begin{align*}
\label{Gint7}
R(z,x) = 2m\,\textrm{sinh}^2t.
\nne
\end{align*}
Then, with $a = a_{\gamma,m/n^2},$ as in \cite{Ku1} (3.23), (the summand of) our Green function integral would become (up to 
a positive constant $C$ depending on normalization of invariant measures)
\begin{align*}
\label{Gint8}
I(\gamma,m,n,v) &= \int_{\Gamma _{a}\setminus \mathbb D}
\int_1^\infty e^ {-2\pi vn^2R(z,a)r}dr/r d\mu(z)\\\nonumber
&= C \mathrm{vol}(\Gamma_{a}\setminus G_{a})
\mathrm{vol}(K)\int_{0}^{\infty}\int_1^\infty e^ 
{-4\pi vm\sinh ^2tr}\sinh t \cosh ^{p-1}t dr/r dt.
\nne
\end{align*}
Hence, for $p=3$ with a similar treatment of the case $m<0,$
we get \textbf{Kudla's formula}
\begin{align*}
\label{Gint9}\nne
I(\gamma,m,v) &= \int_X \Xi(\gamma,m,v) d\mu
=  (1/2) C \sum _{n^2|m}\mathrm{vol}(\Gamma_{a_{m/n^2}}\setminus G_{a_{m/n^2}})\cdot
\mathrm{vol}(K) \cdot I_{\pm}(p,|a|).
\end{align*}
For $p = 3,$ by a routine calculation, we have
the I-Integrals (with $a = 4\pi mv$):
\begin{align*}
\label{epi7a}\nne
 I_{+}(3,a) &= (1/3) \int_{0}^{\infty} e^ {-\alpha w}((w+1)^{3/2}-1)dw/w = (1/3) J_+(3/2,a)\\
I_{-}(3,a) &= (1/3)e^{|a|}\int_{1}^{\infty} e^ {-|a| r}r^{3/2}dr/(r+1) = (1/3)e^{|a|}J_-(3/2,|a|)\\
&= (1/(4|a|^{3/2})\sqrt {\pi}\int_{1}^{\infty} e^ {-|a| r}dr/r^{5/2}.
\end{align*}
The main point in this approach is the formula (\ref{Gint5}) $G = G_{x}A_+K,$ a variant of the Cartan 
decomposition to be found for instance in Theorem 2.4 in Heckman-Schlichtkrull \cite{HS}. 
This is the key to Flensted-Jensen's important integration formula (2.14) in 
his Theorem 2.6 in \cite{FlJ} which Kudla is using above.\\

 The background of Kudla's and Flensted-Jensen's integral formulae  
 is in the following geometric consideration: 
 As to be extracted from \cite{FlJ} Sect.2, we have $G = \textrm{SO}(3.2)$
and for the corresponding Lie algebras, we have
\begin{align*}
\label{Gint10} 
\mathfrak{g} &= \mathfrak{k}+\mathfrak{p}\\
&= \{\left(\begin{smallmatrix}A & \\&B \end{smallmatrix}\right),
A \in M_3(\mathbb{R}) \,\mathrm{skew}, B\in M_2(\mathbb{R})\, \mathrm{skew}\} 
+ \{\left(\begin{smallmatrix} & C\\{}^tC& \end{smallmatrix}\right),
C \in M_{3,2}(\mathbb{R}) \}\\
&= \langle X_{1,2}, X_{1,3}, X_{2,3}, X_{4,5}\rangle 
+ \langle X_{1,4}, X_{1,5}, X_{2,4}, X_{2,5}, X_{3,4}, X_{3,5}\rangle \\
&= \mathfrak{h} + \mathfrak{q}\\
&=  \langle X_{2,3}, X_{4,5}, X_{2,4}, X_{2,5}, X_{3,4}, X_{3,5}\rangle 
+ \langle X_{1,2}, X_{1,3}, X_{1,4}, X_{1,5}\rangle .
\nne
\end{align*}

Here $\mathfrak{k}, \mathfrak{p}$ are the $\pm 1$ eigenspaces 
of the Cartan involution $\tau $ with $\tau X = - {}^tX$  and 
$\mathfrak{h}, \mathfrak{q}$ are the $\pm 1$ eigenspaces 
of the involution $\sigma $ with $\sigma X = E_{1,4}XE_{1,4}.$ 
Corresponding groups are $K = \SO(3)\times \SO(2)$ and 
$H \simeq \SO(2,2).$ We put $L = K\cap H, \mathfrak{l} = \mathrm{Lie}\, L,$ choose $ \mathfrak{b}$ 
maximal abelian in $\mathfrak{p}\cap \mathfrak{q},$ and 
$$
M = Z_L(\mathfrak{b}) = \{\ell \in L; \textrm{Ad}(\ell)B = B \,\,\textrm{for\,all}\,\,B \in  \mathfrak{b} \}.
$$
 In our case, we 
have $\mathfrak{l} =
\mathfrak{k}\cap \mathfrak{h} = \langle X_{2,3}, X_{4,5}\rangle,$ 
 $\mathfrak{p}\cap \mathfrak{q} = \langle X_{1,4}, X_{1,5}\rangle,$ and $\mathfrak{b} = \langle X_{1,4}\rangle.$ Hence $M = \{\mathrm{exp} tX_{2,3};t\in \mathbb{R}\} \simeq \SO(2).$ \\

 In \cite{FlJ} p.261, 
 one observes that the map $L/M\times \mathfrak{b}\rightarrow 
 \mathfrak{p}\cap \mathfrak{q} $ given by 
 \begin{align*}
 \label{Gint11} 
 (lM,B) \mapsto \mathrm{Ad}(l)B
 \nne
 \end{align*}
is a diffeomorphism onto an open dense set. Therefore, 
the map
\begin{align*}
\label{Gint012} 
\Phi:  \mathfrak{p}\cap \mathfrak{h}\times L/M \times  \mathfrak{b} \rightarrow G/K
\nne
\end{align*} 
given by 
$$
\Phi(X,lM,B) = \pi (\mathrm{exp}X \,l\, \mathrm{exp}B), 
$$
where $\pi: G \rightarrow G/K$ is the canonical map, is a diffeomorphism 
unto an open dense set
and also $\Phi': X \mapsto \textrm{exp}X L$ is a diffeomorphism 
of $\mathfrak p \cap \mathfrak h$ unto $H/L.$ \\

The Killing form defines Riemannian (i.e., Euclidean) structures on 
$\mathfrak p \cap \mathfrak h, \mathfrak b^+,$ and $L/M,$ and
one lets the measure on $L/M$ be $\textrm{vol}(L/M)^{-1}$ times the volume element. 
Via Killing form, one has Riemannian structures on $G/K$ and $H/L,$ 
and by their volume elements also measures. Moreover, following Flensted-Jensen, we take measures on $G$ and $H$ such that
\begin{align*}
\label{Flj1} 
\int_G f(x) dx &= \int_{G/K} \int _K f(xk)  dkdxK ,\quad  \int_K dk = 1, \quad \textrm{for}\,f  \in C_c(G) \\
\int_H f(x) dx &= \int_{H/L} \int _L f(xk)  dkdxL ,\quad  \int_L dk = 1, \quad \textrm{for}\, f \in C_c(H) .
\nne
\end{align*}   
Taking the Jacobians $J(X,lM,B) = |\textrm{det}\,d\Phi_{(X,lM,B)}|$ and $J_1(X) =  |\textrm{det}\,d\Phi^1_{(X)}|$
with reference to the respective Riemannian structures, one has for $f  \in C_c(G)$ and $f_1  \in C_c(H)$
\begin{align*}
\label{Flj2} 
\int_{G/K} f(x) dx &= \textrm{vol}\,(L/M)\int_{\mathfrak{p}\cap \mathfrak{h}}\int _{L/M} 
 \int_{\mathfrak b^+}f(\Phi(X,lM,B)) J(X,lM,B)  dBdlM dX \\
 &= \textrm{vol}\,(L/M)\int_{H/M}
 \int_{\mathfrak b^+}f(h \,\textrm{exp}\,B) \delta_1(B)  dB dh 
\nne
\end{align*}   
where $\delta_1(B) =  |\textrm{det}\,d\Phi_{(0,eM,B)}|, \, B \in \mathfrak b^+.$ 
From here (his formula (2.9)), Flensted-Jensen comes 
to the formula (2.14) in his Theorem 2.6
\begin{align*}
\label{Flj3} 
\int_{G} f(g) dg 
 = \textrm{vol}\,(L/M)\int_{K}\int_H
 \int_{\mathfrak b^+}f(k \,\textrm{exp}\,H'h) \delta(H')  dH' dh dk   \quad \textrm{for}\,f  \in C_c(G)
\nne
\end{align*}   
where $\delta$ given by formula (2.12) comes from the $\delta_1.$ 
This also is taken on in a similar way by Kudla-Millson \cite{KMII} (4.35) and (4.37)  as
\begin{align*}
\label{Flj4} 
\int_{\Gamma\backslash G} f(g) dg 
 = \textrm{vol}\,(L/M)\int_{K}\int_{\Gamma\backslash H}
 \int_{\mathfrak b^+}f(h \,\textrm{exp}\,X'k) \delta(X)  dX dh dk   \quad \textrm{for}\,f  \in C_c(G).
\nne
\end{align*}   
In the sequel, we shall use this to interpret Kudla's Green integral.
At first, some remarks concerning the function to be integrated:
 We want to determine $I(\gamma,v,m,n)$ and may use Flensted-Jensen's resp.\,Kudla's formulae, 
 as $R(x,z)$ is left $G_x-$ and right $K-$invariant and depends only on 
 the hyperbolic group $A.$ To do this, we have to assemble and clarify step by step several items. 
 We distinguish between $m>0,$ Case I, and $m<0,$ 
 Case II. Moreover, as above, Case A for $m\in \mathbb{Z}$ and Case B 
 for $m\in \mathbb{Z} + 1/4.$ \\ 
 
Hence, we have the following.\\
 \textbf{Remark.}\label{RemR} 
With $z = g_{x}a_t \cdot z_0,$  one easily has 
\begin{align*}
R(x,z) &= 2m\,\textrm{sinh}^2t \quad \textrm{in Case I with} \quad m>0\\
&= 2|m|\,\textrm{cosh}^2t \quad \textrm{in Case II with} \quad m<0.
\end{align*}

 In Kudla's Green function integrals above for the $(p,2)-$case, the measure $d\mu(z)$ 
on $X$ is 
given by $\Omega^3$ as in \cite{Ku1} (5.8) or \cite{BK} (4.50) with
\begin{align*}
\label{apen67}
\Omega^3 = - \frac{3}{16\pi^3} \,\textrm{det}(y)^{-3}(\frac{i}{2})^3
dz_1\wedge d\bar{z}_1\wedge dz_2\wedge d\bar{z}_2\wedge dz_3\wedge d\bar{z}_3.
\nne
\end{align*}
In the paper \cite{BY} by Bruinier and Yang there is another formula 
relating differential forms and measures. From Proposition 3.4 in \cite{BY}, 
we have 
\begin{align*}
\label{apend67} 
\nne
(d\ell_x)^*(-\Omega)^p = \pm \frac{p!}{(2\pi)^p} \nu_\mathfrak p
\end{align*}
where here $\nu_\mathfrak p$ may seen to be identified with $dx$ in (\ref{Flj2}).
 We take this to realize the Flensted-Jensen formula for an integral of a Green function $f$ as above via
 \begin{align*} \label{fj§}
\int_{\Gamma_x\backslash \mathbb{D}} f(z) d\mu(z) 
 &= (p!/(2\pi)^p)\cdot\textrm{vol}\,(L/M)\int_{K}\int_{\Gamma_x\backslash G_x}
 \int_{\mathfrak b^+}f(h \,\textrm{exp}\,X'k) \delta(X)  dX dh dk,\\
 &= (p!/(2\pi)^p)\cdot\textrm{vol}\,(L/M)\cdot \mathrm{vol}(\Gamma_{x}\setminus G_{x})\cdot 
 I^p_\pm(v,m),
\nne
\end{align*}   
where $d\mu(z) = \Omega^p$ and $\Gamma_x$ and $ I^p_\pm(v,m)$ has to be spezified in each case.\\ 
Now, for $G = \textrm{SO}(3,2),$  
 Case I, and $H = \textrm{SO}(2,2),$ with $a_{0,m/n^2} =: x$ from (\ref{Flj4}), 
 by the usual unfolding we have
 \begin{align*}
\label{repil2}
\nne
I(0,m,v)  
&= (1/2) \int_{\Gamma \backslash \mathbb{D}}  \sum_{x\in L_m} \beta( 2\pi mR(x,z))d\mu(z)\\
&= (1/2)\sum_{x\in \Gamma \backslash L_m}  \int_{\Gamma_x \backslash \mathbb{D}}  \beta( 2\pi mR(x,z))d\mu(z).\\
&= \sum_{n^2|m} I(0,m,n,v),\quad \textrm{with}\\
\cdot I(0,m,n,v) &= (1/2)\int _{\Gamma_{a_{0,m/n^2}} \setminus \mathbb H_2 } \int_1^\infty
e^ {-4\pi vR(na_{0,m/n^2},z)r}dr/rd\mu (z) \\
&=  (1/2)\frac{3!}{(2\pi)^3}
\int_{\Gamma_a\backslash G/K} f(z) dz \\
&=   (1/2)\frac{3!}{(2\pi)^3}\textrm{vol}\,(\textrm{SO}(2))\int_{\Gamma_x\backslash H}
 \int_{0}^\infty\int_1^\infty e^{-4\pi vm (\textrm{sinh}\,t)^2r}|\textrm{sinh}(t)|\textrm{cosh}^2(t)dr/r  dt dh\\
 &=   (1/2)\frac{3!}{(2\pi)^3}\textrm{vol}\,(\textrm{SO}(2))\int_{\Gamma_x\backslash H}dh\,
 I^3_+(v,m)\\
 &=   (3/(4\pi^2))\int_{\Gamma_x\backslash H}dh\,
 I^3_+(v,m),\\
I^3_+(v,m) &= \int_{0}^\infty\int_1^\infty e^{-4\pi vm (\textrm{sinh}\,t)^2r}|\textrm{sinh}(t)|\textrm{cosh}^2(t)dr/r  dt.
\end{align*} 
 
Similarly,  for Case II, $H = \textrm{SO}(3,1),$ 
we get
\begin{align*} 
\label{Apen69} \nne
\cdot I(0,m,n,v) &= (1/2)\int _{\Gamma_{a_{m/n^2}} \setminus \mathbb H_2 } \int_1^\infty
e^ {-4\pi vR(na_{m/n^2},z)r}dr/rd\mu (z) \\
&=  (1/2)\frac{3!}{(2\pi)^3}
\int_{\Gamma_a\backslash G/K} f(z) dz \\
 &=  (1/2) \frac{3!}{(2\pi)^3}\textrm{vol}\,(\textrm{SO}(3)/\textrm{SO}(2))\int_{\Gamma_x\backslash H}dh\,
 I^3_-(v,m).\\
 &=   (3/\pi^2)\int_{\Gamma_x\backslash H}dh\,
 I^3_-(v,m).\\
I^3_-(v,m) &= \int_{0}^\infty\int_1^\infty e^{-4\pi vm (\textrm{cosh}\,t)^2r}(\textrm{sinh}(t))^2\textrm{cosh}(t)dr/r  dt. 
\end{align*} 
As already said, by routine calculations we get
\begin{align*}
 I^3_{+}(v,m) &= (1/3) \int_{0}^{\infty} e^ {-\alpha w}((w+1)^{3/2}-1)dw/w = (1/3) J_+(3/2,a)\\
I^3_{-}(v,m) &= (1/3)e^{-|a|}\int_{1}^{\infty} e^ {-|a| r}r^{3/2}dr/(r+1) = (1/3)e^{-|a|}J_-(3/2,|a|)\\
\end{align*}
Here we used that, e.g.,
from \cite{GHS} (3) and (4), one has $\textrm{vol}\,(\textrm{SO}(2)) = 2\pi,$  $\textrm{vol}\,(\textrm{SO}(3)) = 8\pi^2$ 
and $\textrm{vol}\,(\textrm{SO}(3)/\textrm{SO}(2)) = 4\pi.$
Hence, we have to find the concrete meaning of the measure $dh$ and to determine 
\begin{align*}
\label{apen70}\nne
\textrm{vol}_\ast(\Gamma_x\backslash H) := \int_{\Gamma_x\backslash H}dh.
\end{align*}
As cornerstones we have the two classical results:\\

i) In the book by Elstrodt, Grunewald and Mennicke \cite{EGM}, one has in their Theorem 1.1 in Chapter 7
the following result:\\
Let $K$ be an imaginary quadratic field of discriminant $d_K<0$ and let $\mathcal{O}$ be its ring of 
integers. Then the covolume of the group $\mathbf{PSL}(2,\mathcal{O})$  (in its action 
on the 3-dimensional hyperbolic space $\mathbb{H}^+$) is
\begin{align*}
\label{zuu2}\nne
V_{1,3} := \textrm{vol}( \mathbf{PSL}(2,\mathcal{O})\backslash \mathbb{H}^+) &= \frac{|d_K|^{3/2}}{4\pi ^2}\zeta _K(2)\\
&= \frac{|d_K|^{3/2}}{24}L(2,\chi_K),\, L(s,\chi_K) := \sum _{n>0}(\frac{d_K}{n})n^{-s}.
\end{align*} 
This is also called Humbert's formula and goes back to a result of Humbert from 1919.
Here, with $\mathbb{H}^+ \ni P = (x,y,r),$ the volume is measured with the volume form
\begin{align*}
\label{zuuu2}\nne
dv_{\mathbb{H}^+} = \frac{dxdydr}{r^3}.
\end{align*} 
 
ii) From \cite{HG} p.172 for $m = m_0f^2 > 0, m_0 =d_K$ a fundamental discriminant,  $K = \mathbb{Q}(\sqrt {m_0}),$ we have
\begin{align*}
\label{zuu4}\nne
\textrm{vol}(\textrm{SL}(2,(\mathcal{O}_f,\mathcal{O}^\ast_f)\backslash \mathbb{H}^2) &= 
2f^3 \prod _{p|f}(1- (\frac{m_0}{p})p^{-2}) \zeta _{\mathbb{Q}(\sqrt {m}_0)}(-1)\\
&= 2f^3 \prod _{p|f}(1- (\frac{m_0}{p})p^{-2}) |d_K|^{3/2}2^{-2}\pi ^{-4}\zeta _{\mathbb{Q}(\sqrt {m}_0)}(2)\\
&= f^3 \prod _{p|f}(1- (\frac{m_0}{p})p^{-2}) |d_K|^{3/2}L(2,\chi _{\mathbb{Q}(\sqrt {m}_0)})1/(12\pi^2)
\end{align*}
where, with $\mathbb{H}^2 \ni (z_1 = x_1+iy_1,z_2=x_2+iy_2),$ the volume is measured by
\begin{align*}
\label{zuuu3}\nne
dv_{\textrm{HG}} = \frac{dx_1dy_1dx_2dy_2}{(2\pi y_1y_2)^2},
\end{align*} 
and we used the functional equation
\begin{align*}
\zeta_K(-1) = \zeta_K(2) d_K^{3/2}/(4\pi ^4).
\label{zuu4a}\nne
\end{align*} 
Here $\mathcal{O}_f$ is an order of conductor $f,$ and
\begin{align*}
\label{zuu5}\nne
\textrm{SL}(2,(\mathcal{O}_f,\mathcal{O}^\ast_f)) &= 
\{ \begin{pmatrix}a&b\\c&d \end{pmatrix}; a,d \in \mathcal{O}_f, c \in \mathcal{O}_f^\ast, b \in \mathcal{O}_f^{\ast -1}, ad - bc = 1 \} 
\end{align*} 
with
\begin{align*}
\label{zuu6}\nne
 \mathcal{O}_f^\ast = \{ x \in K; \textrm{Tr}(xa) \in \mathbb{Z} \,\, \textrm{for \, all}\, a \in  \mathcal{O}_f \}.
\end{align*} 
 The formula above goes back to a result by Siegel. As a special case of \cite{S2} (19), one
can conclude that
\begin{align*}
\label{zuu7}\nne
\textrm{vol}(\textrm{SL}(2,\mathcal{O})\backslash \mathbb{H}^2) &= 
(2/\pi^2) |d_K|^{3/2}\zeta _{K}(2) =  |d_K|^{3/2}L(\chi_K,2)/3  =: V_{2,2} \\
&= 8\pi^2\zeta_{K}(-1) \quad\textrm{for}\,d_K > 0
\end{align*}
where, here, the volume is measured by
\begin{align*}
\label{zuuu9}\nne
dv_{\mathbb{H}^2} := \frac{dx_1dy_1dx_2dy_2}{(y_1y_2)^2}.
\end{align*}

In Flensted-Jensen and also in \cite{Ku1}, one usually works with the normalization $\int_K dk = 1.$ 
Moreover, by the geometric meaning of the Eisenstein series (see (\ref{epilo1})) we are 
led to measure the volumes as volumes in the representation space, i.e., here $\mathbb{D}_{22} \simeq \mathbb{H}^2$ resp. $\mathbb{D}_{31}\simeq \mathbb{H}^+.$
We try to put all this together starting by (\ref{Flj4}),
From a rather tedious determination of the unit groups (using Siegel \cite{S3} Section 3 and following the method from 
\cite{EGM}) and the volumes of their fundamental domains, we have 
\begin{align*}
\textrm{vol}_{\textrm{Sie}}(\Gamma^0 \backslash \mathbb{D}_{2,2}) &= 
\textrm{vol}_{\mathbb{H}^2}(\textrm{PSL}_2(\mathcal{O}_f,\mathcal{O}_f^\ast) \backslash \mathbb{H}^2)/4\\
&=  (1/12)|d_F|^{3/2}L(\chi_F,2) \cdot f^3\prod_{p|f} (1 - \chi_{d_F}(p)/p^2),\\ 
\label{SV+}\nne
\textrm{vol}_{\textrm{Sie}}(\Gamma^0 \backslash \mathbb{D}_{1,3}) &=
\textrm{vol}_{\mathbb{H^+}}(\textrm{PSL}_2(\mathcal{O}_f,\mathcal{O}_f^\ast) \backslash \mathbb{H}^+)\\
&=  (1/24)|d_F|^{3/2}L(\chi_F,2) \cdot f^3\prod_{p|f} (1 - \chi_{d_F}(p)/p^2),\
\end{align*}
where $\Gamma^0$ stands for the image of $\Gamma_{a_m}$ 
resp. $\Gamma_{a'_m}$ in the corresponding group $H.$ 
(This comes up to interpret $dh$ in
(\ref{apen70}) as $dv_{\textrm{Sie}}.$)
 To introduce these voluminae into the formulae above, we still have to 
work out the sum of $n^2|d_\gamma^2m$ and the appearance of
$\sigma_{\gamma,m}(5/2)$ in (\ref{SGint2}) .
This goes by the analysis of the volume formula $\mathrm{vol}(\Gamma_{x}\setminus \mathbb{H}^2)$ for the components of the Humbert surfaces in 
\cite{HG}. \\ 

For $x$ with $q(x)>0,$ we have $H = G_x = \SO(2,2)$ and $H = G_x = \SO(3,1)$ for $x$ with $q(x) < 0.$ 
The corresponding homogeneous spaces are $H/K_H = \mathbb{D}_{2,2}$ resp. $\mathbb{D}_{3,1}.$  
Using the $\textrm{SL}_2-$theory and Siegel's method to determine unit groups
in \cite{S3} (on his way to define his \textit{Darstellungsmass}), by some calculation, for $m>0,$ we get

\begin{align*}
\label{epi8}
S_m := \sum _{n^2|m}\mathrm{vol}_{\textrm{BK}}(\Gamma_x\backslash \mathbb{D}_{2,2})
&= \sum _{n^2|m}\mathrm{vol}_{\textrm{BK}}(\Gamma_F\setminus \mathbb{H}^2)f^3\prod_{p|f}(1 - \chi_{D_0}p^{-2}))\\
\nne
\end{align*}
From Theorem 1 (\ref{SGint1a}), we know 
\begin{align*}
\textrm{deg}\,(\mathcal{H}(\gamma,m)) =  - (B/2) C(\gamma,m,0) =
 -  3^{-1} \cdot \pi^{-2} |m|^{3/2}L(\chi_{d_F},2)\sigma_{\gamma,m}(5/2),
\end{align*}
hence $S_m = (B/2) C(\gamma,m,0)$ and for $m>0$ we get
\begin{align*}
\label{repi8}
 \sum _{d|f}(f/d)^3\prod_{p|f/d}(1 - \chi_{D_0}p^{-2})) = f^3\sigma_{\gamma,m}(5/2)\\
\nne
\end{align*}
This formula also follows directly analyzing \cite{BrKu}.
$\hfill\Box$\\

For $m \not= 0,$ In equation (\ref {Apen81}) of the \textbf {Corollary}, we have
\begin{align*}
c'_0(\gamma,m,0,v) = e^{-a/2} ((4/B)\cdot (I(\gamma,m,v) - I^{BK}(\gamma,m,v)) + \ast\, c_0(\gamma,m,0,v) ).
\nne
\end{align*}
This is in line with the result from Ehlen-Sankaran \cite{ES}
\textbf{Theorem 3.6.} For each $z \in \mathbb{D}^0(V),$ the $q-$series
\begin{align*}
- \textrm{log}\,v \varphi_0^\vee + \sum_m (Gr^K_0(m,v) - Gr^B_0(m)) q^m
\label{ES6} 
 \nne
\end{align*} 
is the $q-$expansion of a modular form in $A^!_\kappa(\rho_L^\vee)$ of weight $\kappa = p/2 + 1.$

Rolf Berndt\\
Fachbereich Mathematik, Universit\"at Hamburg, D-20146 Hamburg, Germany\\
berndt@math.uni-hamburg.de
\end{document}